\documentclass[12pt,reqno]{amsart}
\usepackage{amsthm}
\usepackage{amssymb}
\usepackage{graphics}
\usepackage{tikz}
\usepackage{float}
\usetikzlibrary{shapes,backgrounds,calc}
\usepackage{latexsym}
\usepackage{multicol}
\usepackage{verbatim,enumerate}
\usepackage{accents}
\usepackage{cite}
\usepackage{subfig}
\usepackage{multirow}
\usepackage{bigstrut}
\usepackage{amsthm}
\usepackage{amssymb}
\usepackage{graphics}
\usepackage{tikz}
\usetikzlibrary{shapes,backgrounds,calc}
\usepackage{latexsym}
\usepackage{multicol}
\usepackage{verbatim,enumerate}
\usepackage{accents}
\usepackage{cite}
\usepackage{array}
\usepackage[colorlinks=true, linkcolor=blue, citecolor=blue, urlcolor=black]{hyperref}
\usepackage{hyperref}
\usepackage{amsmath, amscd,url}
\usepackage{multirow}
\usepackage{bigstrut}
\usepackage{longtable}

\usepackage{stackengine}

\newcommand{\ind}{\mathop{\mathrm{ind}}}

\newcommand{\F}{\mbox{$\mathbb F$}}	

\usepackage[colorlinks=true, linkcolor=blue, citecolor=blue, urlcolor=black]{hyperref}
\usepackage{hyperref}
\usepackage{amsmath, amscd,url}

\advance\textwidth by 1.3in \advance\oddsidemargin by -.6in \advance\evensidemargin by -.6in
\theoremstyle{definition}

\newtheorem{theorem}{Theorem}[section]

\newtheorem{corollary}[theorem]{Corollary}
\theoremstyle{definition}
\newtheorem{definition}[theorem]{Definition}
\newtheorem{example}[theorem]{Example}

\theoremstyle{remark}

\theoremstyle{definition}

\newcounter{cnt}
 \makeatletter
\def\mydggeometry{\makeatletter\dg@YGRID=1\dg@XGRID=20\unitlength=0.003pt\makeatother}
\makeatother \theoremstyle{remark}

\numberwithin{equation}{section}
\let\bwdg\bigwedge
\def\bigwedge{{\textstyle\bwdg}}

\newcommand{\nc}{\newcommand}
\newcommand{\rnc}{\renewcommand}

\nc{\cal}{\mathcal} \nc{\goth}{\mathfrak} \rnc{\bold}{\mathbf}

\nc\bomega{{\mbox{\boldmath $\omega$}}} \nc\bpsi{{\mbox{\boldmath $\Psi$}}}
 \nc\balpha{{\mbox{\boldmath $\alpha$}}}
 \nc\bpi{{\mbox{\boldmath $\pi$}}}
 \nc\bvpi{{\mbox{\boldmath $\varpi$}}}
\nc\chara{\operatorname{ch}}

  \nc\bxi{{\mbox{\boldmath $\xi$}}}
\nc\bmu{{\mbox{\boldmath $\mu$}}} \nc\bcN{{\mbox{\boldmath $\cal{N}$}}} \nc\bcm{{\mbox{\boldmath $\cal{M}$}}} \nc\blambda{{\mbox{\boldmath
$\lambda$}}}\nc\bnu{{\mbox{\boldmath $\nu$}}}

\makeatletter
\def\section{\def\@secnumfont{\mdseries}\@startsection{section}{1}%
  \z@{.7\linespacing\@plus\linespacing}{.5\linespacing}%
  {\normalfont\scshape\centering}}
\def\subsection{\def\@secnumfont{\bfseries}\@startsection{subsection}{2}%
  {\parindent}{.5\linespacing\@plus.7\linespacing}{-.5em}%
  {\normalfont\bfseries}}
\makeatother

 \nc{\Hom}{\operatorname{Hom}}
  \nc{\mode}{\operatorname{mod}}
\nc{\End}{\operatorname{End}} \nc{\wh}[1]{\widehat{#1}} \nc{\Ext}{\operatorname{Ext}} \nc{\ch}{\text{ch}} \nc{\ev}{\operatorname{ev}}
\nc{\Ob}{\operatorname{Ob}} \nc{\soc}{\operatorname{soc}} \nc{\rad}{\operatorname{rad}} \nc{\head}{\operatorname{head}}

 \nc{\Cal}{\cal} \nc{\Xp}[1]{X^+(#1)} \nc{\Xm}[1]{X^-(#1)}
\nc{\on}{\operatorname} \nc{\Z}{{\bold Z}} \nc{\J}{{\cal J}}  \nc{\Q}{{\bold Q}}

\nc{\N}{{\bold N}}  \nc\boa{\bold a} \nc\bob{\bold b} \nc\boc{\bold c} \nc\bod{\bold d} \nc\boe{\bold e} \nc\bof{\bold f} \nc\bog{\bold g}
\nc\boh{\bold h} \nc\boi{\bold i} \nc\boj{\bold j} \nc\bok{\bold k} \nc\bol{\bold l} \nc\bom{\bold m} \nc\bon{\mathbb n} \nc\boo{\bold o}
\nc\bop{\bold p} \nc\boq{\bold q} \nc\bor{\bold r} \nc\bos{\bold s} \nc\boT{\bold t} \nc\boF{\bold F} \nc\bou{\bold u} \nc\bov{\bold v}
\nc\bow{\bold w} \nc\boz{\bold z}\nc\ba{\bold A} \nc\bb{\bold B} \nc\bc{\mathbb C} \nc\bd{\bold D} \nc\be{\bold E} \nc\bg{\bold
G} \nc\bh{\bold H} \nc\bi{\bold I} \nc\bj{\bold J} \nc\bk{\bold K} \nc\bl{\bold L} \nc\bm{\bold M} \nc\bn{\mathbb N} \nc\bo{\bold O} \nc\bp{\bold
P} \nc\bq{\bold Q} \nc\br{\bold R} \nc\bs{\bold S} \nc\bt{\bold T} \nc\bu{\bold U} \nc\bv{\bold V} \nc\bw{\bold W} \nc\bz{\mathbb Z} \nc\bx{\bold
x} \nc\KR{\bold{KR}} \nc\rk{\bold{rk}} \nc\het{\text{ht }}

\nc\toa{\tilde a} \nc\tob{\tilde b} \nc\toc{\tilde c} \nc\tod{\tilde d} \nc\toe{\tilde e} \nc\tof{\tilde f} \nc\tog{\tilde g} \nc\toh{\tilde h}
\nc\toi{\tilde i} \nc\toj{\tilde j} \nc\tok{\tilde k} \nc\tol{\tilde l} \nc\tom{\tilde m} \nc\ton{\tilde n} \nc\too{\tilde o} \nc\toq{\tilde q}
\nc\tor{\tilde r} \nc\tos{\tilde s} \nc\toT{\tilde t} \nc\tou{\tilde u} \nc\tov{\tilde v} \nc\tow{\tilde w} \nc\toz{\tilde z} \nc\woi{w_{\omega_i}}

\begin{document}
\setcounter{section}{0}
\setcounter{tocdepth}{1}


\title{ON the INDEX divisors OF CERTAIN SEXTIC NUMBER FIELDS}

\author[Anuj Jakhar]{Anuj Jakhar}
\author[Ravi Kalwaniya]{Ravi Kalwaniya}
\address[]{Department of Mathematics, Indian Institute of Technology (IIT) Madras}



\begin{abstract}
\noindent Let $K=\Q(\theta)$ be an algebraic number field with $\theta$ a root of an irreducible quadrinomial $f(x) = x^6+ax^m+bx+c\in\Z[x] $  with $m\in\{2,3,4,5\}$. In the present paper,  we give  some explicit conditions involving only  $a,~b,~c$ and $m$ for which  $K$ is non-monogenic. In each case, we provide  the highest power of a rational prime $p$ dividing index of the field $K$. In particular, we provide a partial answer to the Problem $22$ of
Narkiewicz \cite{Nar} for these number fields. Finally, we illustrate our results through examples.
\end{abstract}
\maketitle
\section{Introduction and statements of results}\label{intro}

Let $f(x)\in\Z[x] $ be a monic irreducible polynomial of degree $n$ over the field  $\Q$ of rational numbers. Let $K=\Q(\theta)$ be an algebraic number field with $\theta$ a root of $f(x)$ and $O_K$ be its ring of algebraic integers.  We denote the index of the subgroup $\Z[\gamma]$ in $O_K$ by $\ind \gamma$. An algebraic number field is called monogenic  if $\ind\alpha=1$ for some $\alpha\in O_K$, i.e., $O_K=\Z[\alpha]$.  In this case $\{1, \alpha,\cdots,\alpha^{n-1}\}$ is an integral basis of $K$; such an integral basis of $K$ is called a power integral basis of $K$. If there does not exist any $\alpha\in O_K$ for which $\ind\alpha=1$, then  $K$ is non-monogenic.  One of the fundamental problems in algebraic number theory is the determination of the monogeneity of an algebraic number field.  The problem of testing the monogeneity of number
fields and constructing power integral bases have been intensively studied (cf.  \cite{ANH1}, \cite{ANH}, \cite{Fun}, \cite{GR}, \cite{MS}, \cite{SS}).  In $1984,$ Funakura \cite{Fun}   gave necessary and sufficient conditions on those integers $m$ for which  the number  field $\Q(m^{1/4})$ is monogenic.
 In $2017$,  Ga\'{a}l and Remete \cite{GR} studied monogenity of  algebraic number fields of the type $\Q(m^{1/n}),$ where  $3\leq n\leq 9$ and $m$ is square free. 
  In \cite{GA}, Ga\'al studied monogenity of number fields defined by some sextic irreducible trinomials. For a rational prime $p$, if the prime ideal factorization of an ideal $pO_K$  is known, then Engstorm \cite{HT} determined the highest power of $p$ dividing the index $i(K)$ of the field $K$ of degree less than or equal to $7.$\\
\indent  In this paper, let $p$ be a rational prime and  $K=\Q(\theta)$ be an algebraic number field  with $\theta$ a root of an irreducible quadrinomial $f(x)=x^6+ax^m+bx+c\in\Z[x]$ with $m\in\{2,3,4,5\},$ then we provide some  sufficient conditions involving only $a,~b,~c$ and $m$ for which $K$ is non-monogenic. Also, we provide the highest power of $p$ dividing the index of the field $K.$ 
  As an application of our results, we provide a family of algebraic number fields which are non-monogenic.\\
  \indent In what follows, for any integer $r$ and a rational prime $p,$ $v_p(r)$ will denote the highest power of $p$ dividing r. Also,  $i(K)$ will stand for the  index of the field $K$ defined by $$ i(K) = \gcd\{\ind \alpha \mid  {\text  {} K=\Q(\alpha)  {\text { and }}      \alpha\in O_K} \}.$$ Note that if for some prime $p,$ $v_p(i(K))\ge 1,$ then the field $K$ is non-monogenic.\\
Precisely, we prove the following result.
\begin{theorem}\label{Th1.2}
Let $K=\Q(\theta)$ be an algebraic number field with $\theta$ a root of an irreducible quadrinomial $f(x)=x^6+ax^m+bx+c\in\Z[x]$ with $m\in\{2,3,4,5\}$. If  $a,b$ and $c+1$ are divisible by $8$ or $9,$ then  $K$ is non-monogenic with $v_2(i(K))=2$ and $v_3(i(K))=1.$  
  \end{theorem}
\vspace{-0.2cm}
 The following corollary follows directly from the above theorem. 
 \begin{corollary}\label{cor1}
   	Let $K=\Q(\theta)$ be an algebraic number field with $\theta$ a root of an irreducible trinomial $f(x)=x^6+ax^m+b\in\Z[x]$ with $m\in\{1,2,3,4,5\}$. If  $a$ and $b+1$ both are divisible by either  $8$ or $9,$  Then $K$ is non-monogenic. \end{corollary}
    \begin{theorem}\label{Th1.3}
     	Let $K=\Q(\theta)$ with $\theta$ a root of an irreducible polynomial $f(x) = x^6+ax^m+bx+c\in \Z[x],$ where $m\in\{2,3,4\}.$ Suppose  $a+(-1)^m, b$ and $c$ are divisible by $8$. Let  $mv_2(b)<(m-1)v_2(c),$ then the following hold:
     	\begin{itemize}
     		\item[(i)] \text{If} $m=2$ and $v_2(a+1-b+c)>3,$ \text{then} $v_2(i(K))=4$
     		\item[(ii)] \text{If} $m=2$ and $v_2(a+1-b+c)=3,$ \text{then} $v_2(i(K))=1$
     		\item[(iii)] \text{If} $m=3$ and $\gcd(v_2(b),2)=1,$ \text{then} $v_2(i(K))=1$
     		\item[(iv)] \text{If} $m=4$ and $\gcd(v_2(b),3)=1,$ \text{then} $v_2(i(K))=2$
     	\end{itemize} 
 In particular, the field $K$ is non-monogenic.\end{theorem}
\begin{theorem}\label{Th1.4}
	Let $K=\Q(\theta)$ with $\theta$ a root of an irreducible polynomial $f(x) = x^6+ax^m+bx+c\in \Z[x],$ where $m\in\{2,3,4\}.$ Suppose $a+(-1)^m, b$ and $c$ are divisible by $9$ and  $mv_3(b)<(m-1)v_3(c)$. Then the following hold:
	\begin{itemize}
     		\item[(i)] \text{If} $m=2,$ \text{then} $v_3(i(K))=1$
     		\item[(ii)] \text{If} $m=3$ and $\gcd(v_3(b),2)=1$, \text{then} $v_3(i(K))=1$
     		\item[(iii)] \text{If} $m=4$ and $\gcd(v_3(b),3)=1$, \text{then} $v_3(i(K))=1$
     	\end{itemize} 
	 In particular, the field $K$ is non-monogenic.
 \end{theorem}
\begin{corollary}\label{cor3}
	Let $K=\Q(\theta)$ with $\theta$ a root of an irreducible polynomial $f(x) = x^6+ax^m+bx+c\in \Z[x],$ where $m\in\{3,4\}.$ Suppose $a+(-1)^m, b$ and $c$ are divisible by either $8$ or $9$. Let $p\in\{2,3\}$ be such that  $mv_p(b)<(m-1)v_p(c),$ and $\gcd(v_p(b),6)=1,$ then by virtue of Theorems \ref{Th1.3} and \ref{Th1.4}, the field $K$ is non-monogenic with $$v_2(i(K))=\begin{cases}
		1& \text{if}~ m=3\\
		2& \text{if}~ m=4\\
	\end{cases}~\text{and}~~v_3(i(K))=1.$$
\end{corollary}
  
 We now provide some examples of non-monogenic number fields.
 \begin{example}
 	Let $K=\Q(\theta)$ with $\theta$ a root of an irreducible polynomial $f(x) = x^6+ax^3+bx+c\in \Z[x].$ Suppose  $a-1, b$ and $c$ are divisible by $8$ and $9.$ If $3v_p(b)<2v_p(c)$ and $v_p(b)$ is odd for $p\in\{2,3\}$, then, by using Theorems \ref{Th1.3} and \ref{Th1.4}, $v_p(i(K))=1.$
 \end{example}
 \begin{example}
 	Let $K=\Q(\theta)$ with $\theta$ a root of a polynomial $f(x) = x^6+ax^2+bx+c\in \Z[x].$ Suppose $a\equiv-7\mod 112$ and $b\equiv 56\mod 112$ and $c\equiv 0\mod 896,$ then, in view of Eisenstein criterion, $f(x)$ is an irreducible polynomial. Hence  by Theorem \ref{Th1.3}, $K$ is non-monogenic and $v_2(i(K))=1$.
 \end{example}
\begin{example}
	Let $K=\Q(\theta)$ with $\theta$ a root of an irreducible polynomial $f(x) = x^6+ax^3+bx+c\in \Z[x].$ Suppose  $a\equiv 1\mod 8,$ $b\equiv 32\mod 64$  and $c\equiv 0\mod 256,$  then, in view of Corollary \ref{cor3}, $v_2(i(K))=1$ and $K$ is non-monogenic.
\end{example}
      
 \section {Preliminary Results}
 Let $K=\Q(\theta)$ be an algebraic number field with $\theta$ a root of a monic irreducible polynomial $f(x)$ belonging to $\Z[x]$. In what follows, $O_K$ will stand for the ring of algebraic integers of $K$.  For a rational prime $p$, let $\F_p$ be the finite field with $p$ elements and $Z_p$ denote the ring of $p$-adic integers. $\overline{f(x)}$ will stand for the polynomial obtained on replacing each coefficient of $f(x)$ modulo $p$.
  
 
 \begin{definition}\label{A} 
The Gauss valuation  of the field $\Q_p(x)$ of rational functions in an indeterminate $x$ which extends the valuation $v_p$ of $\Q_p$ and is defined on $\Q_p[x]$ by \begin{equation*}\label{Gau}
 v_{p,x}(a_0+a_1x+a_2x^2+.....+a_sx^s)= \displaystyle\min\{v_p(a_i),1\leq i\leq s \}, a_i\in \Q_p.
  \end{equation*}
\end{definition}
 \begin{definition}\label{p1.6}
Let $p$ be a rational  prime. Let $\phi(x)\in\Z_p[x]$ be a monic polynomial which is irreducible modulo $p$ and $f(x)\in\Z_p[x]$ be a monic polynomial not divisible by $\phi(x)$. Let $\displaystyle\sum_{i=0}^{n}a_i(x)\phi(x)^i$ with $\deg a_i(x)<\deg\phi(x)$, $a_n(x)\neq 0$ be the $\phi(x)$-expansion of $f(x)$ obtained on dividing it by the successive powers of $\phi(x)$. For each non-zero term $a_i(x)\phi(x)^i$, we associate the point $(n-i,v_{p,x}(a_{i}(x)))$ and form the set $$S=\{(i,~	v_{p,x}(a_{n-i}(x)))\mid 0\leq i\leq n, ~a_{n-i}(x)\neq 0\}.$$ The $\phi$-Newton polygon of $f(x)$ with respect to $p$ is the polygonal path formed by the lower edges along the convex hull of the points of set $S$. The  slopes of the edges form a strictly increasing sequence; these slopes are non-negative as $f(x)$ is a monic polynomial with coefficients in $\Z_p$.
  \end{definition}
   \begin{definition}\label{p1.10}
   Let $\phi(x) \in\Z_p[x]$  be a monic polynomial which is irreducible modulo a rational prime $p$ having a root $\alpha$ in the algebraic closure $\widetilde{\Q}_{p}$ of $\Q_p$. Let $f(x) \in \Z_p[x]$ be a monic polynomial not divisible by $\phi(x)$ with $\phi(x)$-expansion $\phi(x)^n + a_{n-1}(x)\phi(x)^{n-1} + \cdots + a_0(x)$ such that $\overline{f}(x)$ is a power of $\overline{\phi}(x)$. Suppose that the $\phi$-Newton polygon of $f(x)$  consists of a single edge, say $S$, having positive slope denoted by $\frac{l}{e}$ with $l, e$ coprime, i.e.,  $$\min\bigg\{\frac{v_{p,x}(a_{n-i}(x))}{i}~\mid~1\leq i\leq n\bigg\} = \frac{v_{p,x}(a_0(x))}{n} = \frac{l}{e}$$ so that $n$ is divisible by $e$, say $n=et$ and $v_{p,x}(a_{n-ej}(x)) \geq lj$ for $1\leq j\leq t$. Thus the polynomial $b_j(x):=\frac{a_{n-ej}(x)}{p^{lj}}$   has coefficients in $\Z_p$ and hence $b_j(\alpha)\in \Z_p[\alpha]$ for $1\leq j \leq t$.  The polynomial $T(Y)$ in an indeterminate $Y$ defined by   $T(Y) = Y^t + \sum\limits_{j=1}^{t} \overline{b_j}(\overline{\alpha})Y^{t-j}$ having coefficients in $\F_p[\overline{\alpha}]\cong \frac{\F_p[x]}{\langle\phi(x)\rangle}$ is called the residual polynomial of $f
   (x)$ with respect to $(\phi,S)$.
   \end{definition}
   The following definition gives the notion of residual polynomial when $f(x)$ is more general.
   \begin{definition}\label{p1.11} Let $\phi(x), \alpha$ be as in Definition \ref{p1.10}.  Let $g(x)\in \Z_p[x]$ be a monic polynomial not divisible by $\phi(x)$ such that $\overline{g}(x)$ is a power of $\overline{\phi}(x)$. Let $\lambda_1 < \cdots < \lambda_k$ be the slopes of the edges of the $\phi$-Newton polygon of $g(x)$ and $S_i$ denote the edge with slope $\lambda_i$. In view of a classical result proved by Ore (cf. \cite[Theorem 1.5]{CMS}, \cite[Theorem 1.1]{SKS1}), we can write $g(x) = g_1(x)\cdots g_k(x)$, where the $\phi$-Newton polygon of $g_i(x) \in \Z_{{p}}[x]$ has a single edge, say $S_i'$, which is a translate of $S_i$. Let $T_i(Y)$ belonging to ${\F}_{p}[\overline{\alpha}][Y]$ denote the residual polynomial of  $g_i(x)$ with respect to ($\phi,~S_i'$) described as in Definition \ref{p1.10}.  For convenience, the polynomial $T_i(Y)$  will be referred to as the residual  polynomial  of   $g(x)$ with respect to $(\phi,S_i)$. The polynomial $g(x)$ is said to be $p$-regular with respect to $\phi$ if none of the polynomials $T_i(Y)$  has a repeated root in the algebraic closure of $\F_p$, $1\leq i\leq k$. In general, if $F(x)$ belonging to $\Z_p[x]$ is a monic polynomial and $\overline{f}(x) = \overline{\phi}_{1}(x)^{e_1}\cdots\overline{\phi}_r{(x)}^{e_r}$ is its factorization modulo $p$ into irreducible polynomials with each $\phi_i(x)$ belonging to $\Z_p[x]$ monic and $e_i > 0$, then by Hensel's Lemma there exist monic polynomials $f_1(x), \cdots, f_r(x)$ belonging to $\Z_{{p}}[x]$ such that $f(x) = f_1(x)\cdots f_r(x)$ and $\overline{f}_i(x) = \overline{\phi}_i(x)^{e_i}$ for each $i$. The polynomial $f(x)$ is said to be $p$-regular (with respect to $\phi_1, \cdots, \phi_r$) if each $f_i(x)$ is ${p}$-regular with respect to $\phi_i$.
   \end{definition}
   To determine the number of distinct prime ideals of  $O_K$ lying above a rationl prime $p$, we will use  the following theorem which is a weaker version of Theorem 1.2 of \cite{SKS1}.
   \begin{theorem}\label{Th 1}
   Let $L=\Q(\xi)$ be an algebraic number field with $\xi$ satisfying an  irreducible polynomial $g(x)\in \Z[x]$ and $p$ be a rational prime. Let $ \overline{\phi}_{1}(x)^{e_1}\cdots\overline{\phi}_r{(x)}^{e_r}$ be the  factorization of $g(x)$ modulo  $p$ into powers of distinct irreducible polynomials  over $\F_p$ with each $\phi_i(x)\neq g(x)$ belonging to $\Z[x]$  monic. Suppose that  the $\phi_i$-Newton polygon of $g(x)$ has $k_i$ edges, say $S_{ij}$, having slopes $\lambda_{ij}=\frac{l_{ij}}{e_{ij}} $  with $\gcd~(l_{ij},~e_{ij})=1$ for $1\leq j\leq k_i$. If  $T_{ij}(Y) = \prod\limits_{s=1}^{s_{ij}}U_{ijs}(Y)$ is the factorization of  the residual polynomial $T_{ij}(Y)$ into distinct irreducible factors over $\F_p$  with respect to $(\phi_i,~S_{ij})$  for $1\leq j\leq k_i$, then $$pO_L=\displaystyle\prod_{i=1}^{r}\displaystyle\prod_{j=1}^{k_i}\displaystyle\prod_{s=1}^{s_{ij}}\mathfrak p_{ijs}^{e_{ij}},$$ where $\mathfrak p_{ijs}$ are distinct prime ideals of $O_L$ having residual degree $\deg \phi_i(x)\times\deg U_{ijs}(Y).$
   \end{theorem}

Engstrom \cite{HT} calculated the highest power of the rational primes 2 and 3 dividing the index of
any number field K of degree $6$. If $\displaystyle\prod_{i=1}^{k}\mathfrak{p}^{e_i}_i$ is the prime ideal factorization of $pO_K$ with $N(\mathfrak{p}_i)=f_i,$ then the following table  gives the result of these calculations.
\begin{table}[h!]
	\centering
	\begin{tabular}[h!]{m{2cm} m{3cm} m{4cm} m{2cm} m{2cm} }
		\hline
		
Sr. no.&$f_1,f_2,f_3,f_4,f_5$&$e_1,e_2,e_3,e_4,e_5$&$v_2(i(K))$&$v_3(i(K))$\\
		\hline
		$1$&$~2,~2,~1,~1$&$~1,~1,~1,~1$&$~2$& $-$\\
		$2$&$~1,~1,~1,~1$&$~2,~2,~1,~1$&$~2$&$~1$\\
	$3$ &$~1,~1,~1,~1,~1$&$~2,~1,~1,~1,~1$&$~4$&$-$\\
	$4$&$~2,~1,~1,~1$&$~1,~2,~1,~1$&$~1$&$-$\\
	$5$&$~1,~1,~1,~1$&$~3,~1,~1,~1$&$~2$&$~1$\\
	$6$&$~2,~1,~1,~1,~1$&$~1,~1,~1,~1,~1$&$~2$&$~1$\\		
		\hline
	\end{tabular}
\caption{Prime power decomposition of the index of
	sextic number fields.}
\label{Table:1}
\end{table}
\vspace{-0.2in}
    \section{Proof of Theorems \ref{Th1.2}, \ref{Th1.3} and  \ref{Th1.4}}\
    
   \begin{proof}[\textbf{Proof of Theorem \ref{Th1.2}}] \vspace{-0.4cm}
Here  $f(x)\equiv (x+1)^3(x^2+x+1)^3\mod 2$. Set $\phi_1(x)=x^2+x+1$ and $\phi_2(x)=x+1.$ For each $m\in\{2,3,4,5\},$ the $\phi_1(x)$-expansion of $f(x)$ is given by \begin{align*}
m=2:~	& f(x)=\phi_1^3(x)-3x\phi_1^2(x)+(2x+a-2)\phi_1(x)+d_2\\
m=3:~	&f(x)=\phi_1^3(x)-3x\phi_1^2(x)+((a+2)x-a-2)\phi_1(x)+d_3\\
m=4:~	& f(x)=\phi_1^3(x)-(3x+a)\phi_1^2(x)+((2-2a)x-(a+2))\phi_1(x)+d_4\\
m=5:~	&f(x)=\phi_1^3(x)+((a+1)x-2a)\phi_1^2(x)+((a+2)x+3a-2)\phi_1(x)+d_5,\end{align*}
 where \[d_m= 
\begin{cases}
(b-a)x+1+c-a &\text{if}~ m=2\\
bx+1+c+a& \text{if}~ m=3\\
(b+a)x+1+c & \text{if}~ m=4\\
(b+a)x+1+c & \text{if}~ m=5
\end{cases}.
\] As $8$ divides $a, b$ and $c+1$, therefore for each $m\in\{2,3,4,5\}$ we have $v_2(d_m)\ge 3.$ The $\phi_1$-Newton polygon of $f(x)$ is the lower convex hull of the points $(0,~0),$ $(1,~0),$ $(2,~1),$ $(3,~v_2(d_m)).$ The $\phi$-Newton polygon has only two sides of positive slopes. The first side is the line segment joining the points $(1,~0)$ and $(2,~1).$ The second side is the one that joins $(2,~1)$ with $(3,~v_2(d_m)).$ The residual polynomial associated with each side is linear. Next, the $\phi_2$-expansion of $f(x)$ is given by
 \begin{equation}\label{3.1}
 	f(x)=\sum_{i=1}^{6}	\left((-1)^i{6\choose i}+a(-1)^{m-i}{m \choose i}\right)\phi_2^i(x)+(a(-1)^m-b+1+c).
 \end{equation} The $\phi_2$-Newton polygon of $f(x)$  has two sides of positive slopes. The first side is the line segment joining the points $(4,~0)$ and $(5,~1)$. The second side is the one that joins $(5,~1)$ with $(6,~v_2(a(-1)^m-b+1+c)).$ The residual polynomial attached to each side is linear. Thus Theorem \ref{Th 1} implies that $\phi_1$ provides two prime ideals, say $\mathfrak{p}_1$ and $\mathfrak{p}_2$, of residual degree $2$ each and $\phi_2(x)$ give two prime ideals, say $\mathfrak{p}_3$ and $\mathfrak{p}_4$, of residual degree $1$ each. So $2O_K=\mathfrak{p}_1\mathfrak{p}_2\mathfrak{p}_3\mathfrak{p}_4.$  Hence by using Table \ref{Table:1}, we have $v_2(i(K))=2.$ \\ \indent  When $9$ divides $a,b$ and $c+1,$ then $f(x)\equiv (x+1)^3(x-1)^3\mod 3.$ Take $\phi_3(x)=x-1.$ The  $\phi_3(x)$-expansion is given by
\begin{equation}\label{3.2}
	f(x)=\sum_{i=1}^{6}	\left({6\choose i}+a{m \choose i}\right)\phi_3^i(x)+(a+b+c+1).
\end{equation}  Let \[d_j= 
\begin{cases}
a(-1)^m-b+c+1 &\text{if}~ j=2\\
a+b+c+1& \text{if}~ j=3
\end{cases}.
\] Then keeping in mind the  $\phi_2(x)$-expansion given in \eqref{3.1} and $\phi_3(x)$-expansion of $f(x)$, we observe that for each $j=2,3,$ the $\phi_j$-Newton polygon of $f(x)$ being the lower convex hull  of the points $(0,~0),$ $(1,~1),$ $(2,~1),$ $(3,~0),$ $(4,~1),$ $(5,~1)$ and $(6,~v_3(d_j))$ has two sides of positive slopes. The residual polynomial attached to each side is linear. Therefore in view of Theorem \ref{Th 1}, $3O_K=\mathfrak{p}^2_1\mathfrak{p}^2_2\mathfrak{p}_3\mathfrak{p}_4,$ where residual degree of $\mathfrak{p}_i$ for $i=1,2,3,4,$ is $1.$  Thus in view of the Table \ref{Table:1}, $v_3(i(K))=1.$ Hence the field $K$ is non-monogenic. This completes the proof.
\end{proof}

\begin{proof}[\textbf{Proof of Theorem \ref{Th1.3}}] Here  \[f(x)\equiv 
	\begin{cases}
		x^2(x+1)^4 \mod 2       & \text{if}~ m=2\\
		x^3(x+1)(x^2+x+1)\mod 2 &\text{if}~ m=3\\
		x^4(x+1)^2 \mod 2       & \text{if} ~ m=4
	\end{cases}
	\] The $x$-Newton polygon of $f(x)$ is the lower convex hull of the points $(0,~0),$ $(6-m,~0),$ $(5,~v_2(b))$  and $(6,~v_2(c)).$ As $mv_2(b)<(m-1)v_2(c),$ therefore for each $m\in\{2,3,4\},$ it has two sides of positive slopes. The first side is the line segment joining the points $(6-m,~0)$ and $(5,~v_2(b))$ and the second side joins $(5,~v_2(b))$ to $(6,~v_2(c)).$ Keeping in mind the  hypothesis ($\gcd(v_2(b),2)=1$ for $m=3$; $\gcd(v_2(b),3)=1$ for $m=4$), we observe that the residual polynomial corresponding to each side is linear. By using  $x+1$ expansion given in \eqref{3.1}, if  $m=2$, then we see  that the $x+1$-Newton polygon of $f(x)$ has either two or three sides (depending on $v_2(a+1-b+c)$) of positive slopes  and if $m=3,$ then it has two sides of positive slopes. The residual polynomial attached to each side is linear.  If $g(x)\in\{x+1,x^2+x+1\},$  then for $m=3,$ $g(x)$-Newton polygon has a single side of positive slope and the residual polynomial attached to this side is linear. Therefore by virtue of Theorem \ref{Th 1}, and Table \ref{Table:1}, we have
	\begin{center}
		\begin{longtable}[h!]{|m{2cm}|m{4cm}| m{4cm}|m{2cm}|}
			\hline
			Case &$2O_K$&$f_1,f_2,f_3,f_4,f_5$ & $v_2(i(K))$\\
			\hline
			$m=2$&$\mathfrak{p}^2_1\mathfrak{p}_2\mathfrak{p}_3\mathfrak{p}_4\mathfrak{p}_5,$  &$~1,~1,~1,~1,~1$&$~4$\\
			&or&&\\
			& $\mathfrak{p}_1\mathfrak{p}^2_2\mathfrak{p}_3\mathfrak{p}_4,$&$~2,~1,~1,~1$&$~1$\\
			$m=3$&$\mathfrak{p}_1\mathfrak{p}^2_2\mathfrak{p}_3\mathfrak{p}_4$& $~2,~1,~1,~1$&$~1$\\
			$m=4$&$\mathfrak{p}^3_1\mathfrak{p}_2\mathfrak{p}_3\mathfrak{p}_4$&$~1,~1,~1,~1$&$~2$\\		
			\hline
		\end{longtable}
	\end{center}\vspace{-1cm}
	  where $f_i$ denote the residual degree of $\mathfrak{p}_i.$
\end{proof}
	 \begin{proof}[\textbf{Proof of Theorem \ref{Th1.4}}] Clearly \[f(x)\equiv 
  	\begin{cases}
  		x^2(x-1)(x+1)(x^2+1) \mod 3       & \text{if}~ m=2\\
  		x^3(x+1)^3           \mod 3 &\text{if}~ m=3\\
  		x^4(x-1)(x+1)        \mod 3       & \text{if} ~ m=4
  	\end{cases}
  	\] The $x$-Newton polygon of $f(x)$ is the lower convex hull of the points $(0,~0),$ $(6-m,~0),$ $(5,~v_3(b))$  and $(6,~v_3(c)).$ As $mv_3(b)<(m-1)v_3(c),$ therefore the $x$-Newton polygon of $f(x)$ has two sides of positive slopes. The first side is the line segment joining the points $(6-m,~0)$ and $(5,~v_3(b))$ and the second side  joins $(5,~v_3(b))$ to $(6,~v_3(c)).$ Using hypothesis, we see that the residual polynomial corresponding to each side is linear. For $m=3,$ the $x+1$-Newton polygon of $f(x)$ has two sides of positive slopes and each side has linear residual polynomial. If $h(x)\neq x$ is an irreducible factor of $f(x)$ in modulo $3,$ then for $m\in\{2,4\},$ $h(x)$-Newton polygon has a single side of positive slope and the residual polynomial attached to this side is linear. In view of Theorem \ref{Th 1} and Table \ref{Table:1}, we have 
  	\begin{center}
  		\begin{longtable}[h!]{|m{2cm}|m{4cm}| m{4cm}|m{2cm}|}
  			\hline
  			Case &$3O_K$&$f_1,f_2,f_3,f_4,f_5$ & $v_3(i(K))$\\
  			\hline
  			$m=2$&$\mathfrak{p}_1\mathfrak{p}_2\mathfrak{p}_3\mathfrak{p}_4\mathfrak{p}_5,$&$~2,~1,~1,~1,~1$&$~1$\\
  			$m=3$&$\mathfrak{p}^2_1\mathfrak{p}^2_2\mathfrak{p}_3\mathfrak{p}_4$& $~1,~1,~1,~1$&$~1$\\
  			$m=4$&$\mathfrak{p}^3_1\mathfrak{p}_2\mathfrak{p}_3\mathfrak{p}_4$&$~1,~1,~1,~1$&$~1$\\		
  			\hline
  		\end{longtable}
  	\end{center}\vspace{-1cm}
  	 Hence the field $K$ is non-monogenic.
  \end{proof}

 \newpage


\end{document}